\documentclass[11pt,a4paper]{article}
\usepackage[utf8]{inputenc}
\usepackage[spanish]{babel}
\usepackage{pst-grad} 
\usepackage{pst-plot} 
\usepackage{pstricks}
\usepackage{amsmath,amssymb,mathrsfs,amsthm}
\usepackage[dvips]{graphicx}
\usepackage{epsfig}

\newcommand{\RR}{\mathbb R}

\newcommand{\TT}{\mathbb T}
\newcommand{\pat}{\partial_t}
\newcommand{\pax}{\partial_x}

\newcommand{\reps}{\rho_\epsilon}

\usepackage[pdfauthor={Rafael Granero Belinch\'on y Jos\'e Manuel Moreno Valderrama},pdftitle={La ecuaci\'on de Burgers como un paso previo al estudio de los fluidos incompresibles},bookmarks,colorlinks]{hyperref}
\hypersetup{colorlinks=true}

\newcounter{comentcount}
\newcounter{teocount}
\newtheorem{lem}{Lema}
\newtheorem{prop}{Proposici\'on}
\newtheorem{teo}[teocount]{Teorema}

\newenvironment{coment}
{\stepcounter{comentcount} {\bf\tt Comentario} {\bf\tt\arabic{comentcount}} }{ }
\title{La ecuaci\'on de Burgers como un paso previo al estudio de los fluidos incompresibles}
\author{Rafael Granero Belinch\'on$^{\mbox{{\footnotesize 1}}}$ y Jos\'e Manuel Moreno Valderrama$^{\mbox{{\footnotesize 2}}}$}

\begin{document}

\maketitle 
\footnotetext[1]{Email: \texttt{r.granero@icmat.es}\\
Consejo Superior de Investigaciones Cient\'ificas\\
Instituto de Ciencias Matem\'aticas (CSIC-UAM-UC3M-UCM)\\
C/Nicol\'as Cabrera, 13-15,\\
Campus de Cantoblanco,\\
28049 - Madrid}
\footnotetext[2]{Email: \texttt{josemanuel.moreno@estudiante.uam.es}\\
Universidad Aut\'onoma de Madrid\\
Campus de Cantoblanco,\\
28049 - Madrid}
\vspace{0.3cm}

\vspace{0.3cm}

\begin{abstract}
En este art\'iculo se presentan las ecuaciones de Euler y de Navier-Stokes, las m\'as b\'asicas de las ecuaciones de los fluidos incompresibles, así como unos modelos simplificados de dichos problemas como pueden ser la ecuaci\'on quasigeostr\'ofica u otros escalares activos. Las t\'ecnicas matem\'aticas para obtener que estas ecuaciones est\'an bien puestas en el sentido de Hadamard y otras propiedades cualitativas (Principios del M\'aximo, formaci\'on de singularidades...) se ilustran en el caso sensiblemente m\'as sencillo de la ecuaci\'on de Burgers con disipaci\'on no-local. Se adjunta una secci\'on dedicada a los m\'etodos num\'ericos usados para aproximar soluciones a estos problemas. Este trabajo tiene su origen en una serie de clases que impartimos en la escuela JAE-Intro del CSIC durante el verano del curso 2010-2011 y por lo tanto se centra en la popularizaci\'on y divulgaci\'on de la f\'isica involucrada y en la explicaci\'on detallada de las ideas m\'as abstractas en los argumentos puramente matem\'aticos.
\end{abstract}

\vspace{0.3cm}


\textbf{Palabras clave}: Ecuaciones de Euler, ecuaciones de Navier-Stokes, ecuaci\'on quasi-geostr\'ofica, medios porosos, Ley de Darcy, interfase.

\textbf{Agradecimientos}: R.Granero est\'a financiado por el proyecto MTM2008-03754 del Ministerio de Ciencia e Innovaci\'on (MICINN).

\tableofcontents

\section{Introducci\'on y motivaci\'on}
El estudio de las ecuaciones de los fluidos incompresibles tiene cada vez un mayor inter\'es, tanto desde el punto de vista m\'as te\'orico (integrales singulares...) como desde el enfoque m\'as aplicado (simulaciones num\'ericas...). 

Las ecuaciones que aparecen modelizando problemas de mec\'anica de fluidos son variadas, pero las m\'as importantes son las de Euler y Navier-Stokes. De hecho por demostrar (o refutar) la existencia global de soluci\'on cl\'asica para Navier-Stokes el Instituto Clay otorga un premio de un mill\'on de d\'olares. 

Consideramos como dominio espacial el plano o el espacio enteros, $i.e.\;\RR^d$ con $d=2,3$ y exigimos que la velocidad del fluido sea cero en el infinito. Sea adem\'as $f$ el campo de velocidades inicial. Entonces las ecuaciones de Euler (1707-1783) para la velocidad $u$ de un fluido incompresible son
\begin{equation}
\left\{\begin{array}{l}
\stackrel{\text{Masa}}{\overbrace{\rho}}\stackrel{\text{Aceleraci\'on}}{\overbrace{(\partial_t u +(u\cdot \nabla) u)}}=\stackrel{\text{Fuerzas internas}}{\overbrace{-\nabla p}}+\stackrel{\text{Fuerzas externas}}{\overbrace{F}}\;\; (\text{Conservaci\'on del momento}),\\
\nabla\cdot u=0,\;\; (\text{Conservaci\'on de la masa}) .
\end{array}\right.
\label{Euler}
\end{equation}
con $\nabla=(\partial_{x_1},...,\partial_{x_d})$ y $u=(u_1,...,u_d)$. Este sistema de ecuaciones es la segunda ley de Newton en el caso de un continuo de part\'iculas. 

Las ecuaciones de Navier (1785-1836) y Stokes (1819-1903) para un fluido incompresible con densidad $\rho$ son
\begin{equation}
\left\{\begin{array}{l}
\stackrel{\text{Masa}}{\overbrace{\rho}}\stackrel{\text{Aceleraci\'on}}{\overbrace{(\partial_t u +(u\cdot \nabla) u)}}=\stackrel{\text{Fuerzas internas}}{\overbrace{-\nabla p+\nu \Delta u}}+\stackrel{\text{Fuerzas externas}}{\overbrace{F}}\;\; (\text{Conservaci\'on del momento}),\\
\nabla\cdot u=0,\;\; (\text{Conservaci\'on de la masa}),
\end{array}\right.
\label{NS}
\end{equation}
donde $\Delta=\sum_{i=1}^d\partial_{x_i}^2$.

En las ecuaciones de Navier-Stokes se ha a\~ nadido el rozamiento entre part\'iculas del fluido modeliz\'andolo con un laplaciano. 

Observamos que es un sistema de evoluci\'on \emph{'extra\~ no'} en el sentido de que la derivada temporal de la presi\'on $p$ no aparece. Eso nos indica que la presi\'on se puede obtener de la velocidad $u$. Para ver esto basta tomar la divergencia de la ecuaci\'on de conservaci\'on del momento,
$$
\nabla\cdot((u\cdot\nabla)u)=-\Delta p.
$$
Ahora podemos utilizar la funci\'on de Green para el Laplaciano (que para $\RR^d$ es conocida) y obtener $p=G(u)$. Tomamos el gradiente en esta expresi\'on y obtenemos un sistema de ecuaciones no-locales cerrado para $u$ (donde ahora nos restringimos a velocidades incompresibles). Adem\'as hemos obtenido que la presi\'on act\'ua como un multiplicador de Lagrange para la restricci\'on $\nabla\cdot u=0$.
Los resultados de los que se dispone cierran la teor\'ia en $d=2$. En el plano se conoce la existencia global (en tiempo) de soluciones cl\'asicas. Para el problema completo $d=3$ se tiene un teorema de existencia local (en tiempo) de soluci\'on cl\'asica (ver \cite{L}). Esta diferencia es porque las ecuaciones en los casos bidimensional y tridimensional son sensiblemente distintas. Para verlo necesitamos definir una cantidad que nos de informaci\'on sobre cu\'anto \emph{'gira'} el fluido. Esta cantidad que denotaremos $\omega$ se llama vorticidad y es el rotacional de la velocidad, $i.e.$
$$
\text{rot }u=\omega,
$$
con $u$ el campo de velocidades del fluido. En concreto se tiene que las ecuaciones para la vorticidad son distintas: para un fluido no viscoso e incompresible\footnote{Estos fluidos se conocen como \emph{fluidos ideales}.}, es decir, que sigue las ecuaciones de Euler, en el caso $d=2$ se tiene
\begin{equation}
\pat \omega+ u\cdot\nabla\omega=0,
\label{vort2} 
\end{equation}
que es una ecuaci\'on de transporte, mientras que en el caso $d=3$ se tienen las siguientes ecuaciones
\begin{equation}
\pat \omega+ (u\cdot\nabla)\omega=(\omega\cdot \nabla)u.
\label{vort3} 
\end{equation}
En el caso de considerar un fluido viscoso se ha de a\~ nadir un t\'ermino $\Delta \omega$.

\begin{figure}[h]
	\centering
		\includegraphics[scale=0.3]{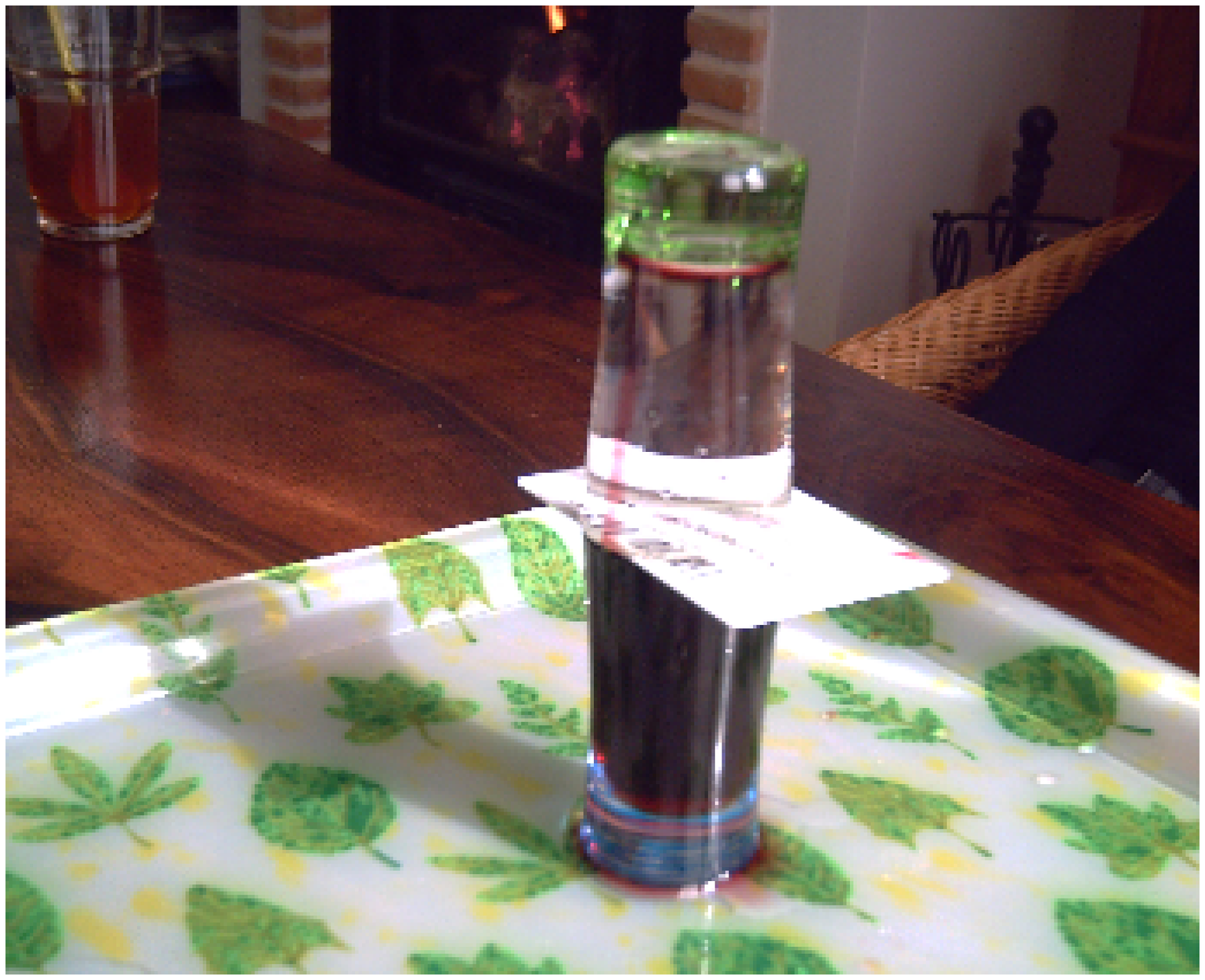} \includegraphics[scale=0.3]{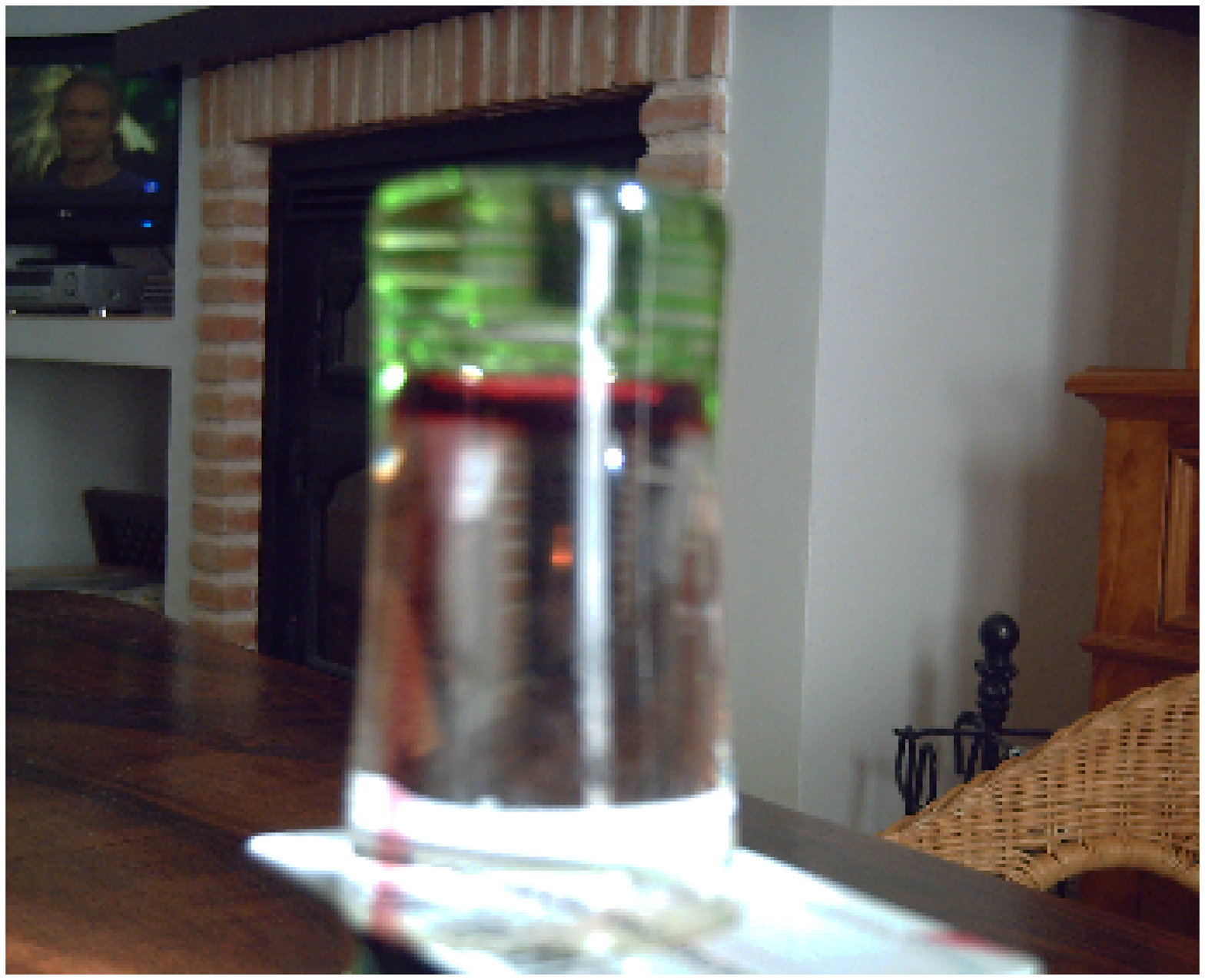}
		\caption{Experimento.}
\label{figura1}
\end{figure}

La diferencia entre ambos casos aparece tambi\'en en las longitudes de onda que est\'an \'intimamente relacionadas con el fen\'omeno de la \emph{'turbulencia'} (ver \cite{DG}). La turbulencia tiene como efecto principal facilitar que dos fluidos se mezclen, por lo tanto, llegados a este punto podemos \emph{'experimentar'} un teorema. Para este peque\~ no juego necesitamos dos vasos peque\~ nos llenos hasta arriba uno de ellos de agua y el otro de vino. La cuesti\'on es: \textquestiondown c\'omo conseguimos cambiar los l\'iquidos de vaso sin usar un tercer recipiente y sin que se mezclen? Para responder a esta pregunta hemos de conocer c\'omo es la turbulencia en tres dimensiones y qu\'e diferencia hay con dos dimensiones. As\'i, si conseguimos una manera de \emph{reducir el problema tridimensional a uno bidimensional} hemos acabado, porque en dos dimensiones \emph{'no hay turbulencia'} y \'esta es la culpable de que los l\'iquidos se mezclen. Para conseguir esta reducci\'on en las dimensiones lo que hacemos es tapar el vaso de agua con un carn\'e y colocarlo con cuidado encima del vaso de vino. Si lo hemos hecho bien no se ha salido ni una gota. Ahora abrimos una rendija min\'uscula entre los vasos y el carn\'e. El agua es m\'as densa, por lo tanto comenzar\'a a bajar mientras que el vino subir\'a...\textexclamdown y todo esto sin mezclarse! (ver Figura \ref{figura1})

Tras esta excursi\'on por las ciencias experimentales volvamos a las matem\'aticas. La vorticidad $\omega$ es una cantidad que aparece en el conocido criterio de existencia global de Beale-Kato-Majda (ver \cite{bkm}).
\begin{teo}[Beale-Kato-Majda]
Sean $T^*$ y $M$ dos constantes tales que 
$$
\int_0^T||\omega||_{L^\infty}\leq M\;\;\forall T<T^*,
$$
entonces una soluci\'on cl\'asica de las ecuaciones de Euler \eqref{Euler} existe al menos hasta tiempo $T^*$. Adem\'as, si $T^{max}$ es el tiempo m\'aximo de existencia (es decir, aparece una singularidad), entonces
$$
\lim_{T\rightarrow T^{max}}\int_0^T||\omega||_{L^\infty}\rightarrow\infty.
$$
\label{teo1}
\end{teo}
Este teorema nos dice que si controlamos \emph{'lo que gira'} el fluido entonces no hay singularidades, por lo tanto conocer c\'omo se comporta la vorticidad es crucial para intentar entender qu\'e hace el fluido. En efecto, tambi\'en es interesante porque si la conocemos podemos recuperar la velocidad gracias a la f\'ormula de Biot-Savart (ver \cite{bertozzi-Majda}). En el caso bidimensional la ley de Biot-Savart es
\begin{equation}
u=\int_{\RR^2}K(x-y)\omega(y)dy=K(\omega),
\label{BS}
\end{equation}
con 
$$
K(x)=\frac{1}{2\pi}\left(\frac{-x_2}{|x|^2},\frac{x_1}{|x|^2}\right).
$$
Ahora bien, con lo que acabamos de mencionar podemos ver que el problema de las ecuaciones de Euler \eqref{Euler} en el caso bidimensional podemos formularlo de manera equivalente como un \emph{'escalar activo'}, es decir, un escalar que es transportado por el fluido de manera que adem\'as podemos recuperar la velocidad del fluido si conocemos el escalar. En efecto: recordemos que ten\'iamos la ecuaci\'on para la vorticidad (que en dos dimensiones es un escalar) \eqref{vort2}, esto unido a $u=K(\omega)$ cierra el problema para la vorticidad. Adem\'as si asumimos que nuestro dato inicial est\'a acotado, usando el Teorema \ref{teo1}, tenemos que la soluci\'on cl\'asica del sistema \eqref{Euler} existe globalmente.

El estudio matem\'atico de los escalares activos tiene una gran relevancia en cuanto que son sistemas \emph{sencillos} que conservan el car\'acter no local de un fluido incompresible. Adem\'as de que hay multitud de aplicaciones donde el problema f\'isico se puede modelizar con un escalar activo. Por ejemplo tenemos el caso de la ecuaci\'on \emph{Quasigeostr\'ofica} (ver \cite{cv}, \cite{cor2}, \cite{G}, \cite{ccw} y las referencias all\'i expuestas) o de la \emph{Ley de Darcy} (ver \cite{c-g07},\cite{c-g09}, \cite{c-g-o08}, \cite{c-c-g10}, \cite{bear}, \cite{M}). La ecuaci\'on quasigeostr\'ofica modeliza la evoluci\'on de la temperatura de grandes masas de aire en grandes escalas espaciales y es muy estudiada como modelo de la \emph{'frontog\'enesis'} (la formaci\'on de frentes de aire a distinta temperatura). Este problema es de inter\'es en meteorolog\'ia, porque ya se sabe que 
\newpage
\begin{quote}                                                                                                                                                                                                                                                                                                                                                                                                                                                                                                                                                                                                                                  
La falta de acierto de quienes predicen el tiempo se ha hecho ya proverbial, y sin embargo no hay ning\'un meteor\'ologo competente que no opine que los procesos atmosf\'ericos est\'an causalmente determinados. 
\begin{flushright}
Max Planck (extra\'ido de \cite{P})
\end{flushright}
\end{quote}

\noindent
La ecuaci\'on quasigeostr\'ofica en dos dimensiones espaciales es 
\begin{equation}
\left\{\begin{array}{l}
\pat \theta+ u\cdot\nabla\theta=0,\\
u=R^\perp u=(-R_2\theta,R_1\theta),
\end{array}\right.
\label{QG}
\end{equation}
donde $R_i$ es la transformada de Riesz $i-$\'esima (ver \cite{St}) y $\theta$ es la temperatura del aire. 

La Ley de Darcy modeliza un fluido incompresible que se mueve a bajas velocidades por un medio poroso. As\'i si $\rho(x,t)$ es la densidad del fluido se tiene el problema
\begin{equation*}
\left\{\begin{array}{l}
\pat \rho+ u\cdot\nabla\rho=0,\\
u=-(\nabla p+ge_2\rho),\\
\nabla\cdot u=0.
\end{array}\right.
\end{equation*}
Este problema puede reducirse a un escalar activo tomando el rotacional dos veces en la ecuaci\'on, obteniendo 
\begin{equation}
\left\{\begin{array}{l}
\pat \rho+ u\cdot\nabla\rho=0,\\
u=(-\Delta)^{-1}( \text{rot rot }ge_2\rho).\\
\end{array}\right.
\label{Darcy}
\end{equation}
Podemos modelizar la transferencia de calor interno del fluido en un medio poroso con la misma ecuaci\'on \eqref{Darcy} si a $\rho$ le damos el sentido de una temperatura (ver \cite{bn}). Por lo tanto, tanto en el caso de la ecuaci\'on \eqref{QG} y de \eqref{Darcy} puede interesarnos a\~ nadir un t\'ermino de difusi\'on del calor. Sin embargo estos t\'erminos de difusi\'on no tienen por qu\'e ser el t\'ipico laplaciano, puede ser necesario a\~ nadir una \emph{'potencia fraccionaria del laplaciano'}. El laplaciano en el espacio de Fourier (con variables $\xi$) tiene una expresi\'on sencilla, es un multiplicador:
$$
\widehat{-\Delta u}=|\xi|^2\hat{u}.
$$
Podemos definir el operador $\Lambda=\sqrt{-\Delta}$ de la siguiente manera
\begin{equation}
\widehat{\Lambda u}=|\xi|\hat{u},
\label{Lambda} 
\end{equation}
y equivalentemente
\begin{equation}
\widehat{\Lambda^\alpha u}=|\xi|^\alpha\hat{u}.
\label{Lambdaalpha} 
\end{equation}
Notemos que tambi\'en podemos escribir el resultado de aplicar el laplaciano fraccionario como la siguiente convoluci\'on:
\begin{equation}
\Lambda^{\alpha} u(x) = \beta({\alpha},d)\text{P.V.}\int_{\mathbb{R}^{d}}\frac{u(x)-u(y)}{|x-y|^{d+\alpha}}dy
\label{conv}
\end{equation}
donde $\beta(\alpha,d)$ es una constante de normalizaci\'on. Los operadores \eqref{Lambda} y \eqref{Lambdaalpha} (o su versi\'on \eqref{conv}) son lo que nosotros entendemos por \emph{'potencias fraccionarias del laplaciano'}. As\'i las ecuaciones \eqref{QG} y \eqref{Darcy} con difusi\'on no local son
\begin{equation}
\left\{\begin{array}{l}
\pat \theta+ u\cdot\nabla\theta=-\gamma \Lambda^\alpha \theta,\\
u=R^\perp u=(-R_2\theta,R_1\theta),
\end{array}\right.
\label{QGvis}
\end{equation}
y
\begin{equation}
\left\{\begin{array}{l}
\pat \rho+ u\cdot\nabla\rho=-\gamma \Lambda^\alpha \rho,\\
u=(R_1R_2, -R_1^2)\rho.\\
\end{array}\right.
\label{Darcyvis}
\end{equation}

Para irnos aproximando a estos problemas podemos plantearnos otros m\'as sencillos o simplificados. Por ejemplo, consideremos un escalar activo
\begin{equation}
\left\{\begin{array}{l}
\pat \eta+ u\cdot\nabla\eta=-\gamma \Lambda^\alpha \eta,\\
u=T(\eta),\\
\end{array}\right.
\label{EA}
\end{equation}
donde $T$ es un operador integral singular que adem\'as nos garantiza que $\nabla\cdot u=0$. Queremos simplificar este problema de manera que sea f\'acilmente abordable, pero no tiene que estar tan simplificado que no nos de ninguna informaci\'on. Lo primero que hacemos es reducir el n\'umero de variables espaciales a una, por lo tanto $\nabla=\pax$. Adem\'as podemos simplificar $T$ tom\'andolo igual a la identidad, es decir, pierde su car\'acter integral singular. La ecuaci\'on resultante de estas simplificaciones es la \emph{ecuaci\'on de Burgers viscosa} si $\gamma>0$ y la \emph{ecuaci\'on de Burgers no viscosa} si $\gamma=0$ (ver \cite{KNS} y \cite{DDL}):
\begin{equation}
\pat \eta+\eta\pax\eta=-\gamma\Lambda^\alpha \eta.
\label{Burgers} 
\end{equation}
Como la condici\'on $\nabla\cdot u=0$ en una dimensi\'on no tiene sentido f\'isico (en una dimensi\'on no hay choques) con esta ecuaci\'on tenemos un modelo unidimensional de las ecuaciones \eqref{Euler} y \eqref{NS}. As\'i nuestra ecuaci\'on \eqref{Burgers} es el modelo m\'as sencillo que nos da informaci\'on tanto sobre las ecuaciones \eqref{Euler}, \eqref{NS} como \eqref{EA}.

Este texto est\'a organizado de la siguiente manera: en la secci\'on 2 probaremos la existencia local de soluci\'on cl\'asica para toda la familia de ecuaciones \eqref{Burgers}. En la secci\'on 3 obtendremos una ley de conservaci\'on y unos principios del m\'aximo que se tienen para las soluciones cl\'asicas de la ecuaci\'on \eqref{Burgers}. En la secci\'on 4 daremos un criterio de existencia de soluci\'on cl\'asica an\'alogo al de Beale-Kato-Majda. En la secci\'on 5 veremos que hay \emph{blow up} para $\gamma=0$ en \eqref{Burgers} y en la secci\'on 6 haremos simulaciones num\'ericas de la soluci\'on de \eqref{Burgers}.

\section{Existencia local de la soluci\'on cl\'asica}
El problema que estudiaremos en esta secci\'on y las siguientes es la ecuaci\'on \eqref{Burgers} (donde cambiamos la notaci\'on $\eta$ por la m\'as corriente $u$):
\begin{equation}
\left\{\begin{array}{l}
\pat u+ u\pax u=-\gamma\Lambda^\alpha u,\;\; (t,x)\in[0,T]\times\RR\\
u(0,x)=f(x),\\
\end{array}\right.
\label{Burg}
\end{equation}
donde $T>0$ es el tiempo de existencia, $0<\alpha\leq 2$ y $\gamma\geq0$. Sobre el dato inicial impondremos las condiciones necesarias cuando veamos el teorema de existencia local, pero por el momento podemos suponer que $f\in C^2$. Nos restringiremos a las \emph{soluciones cl\'asicas} o \emph{soluciones fuertes}. Este tipo de soluciones es tan suave como sea necesario para dar el sentido usual a las derivadas parciales que aparecen en la ecuaci\'on. As\'i para la ecuaci\'on de Burgers queremos que $u$ tenga al menos una derivada en $t$ y dos en $x$. Es decir, queremos que $u\in C^{2,1}_{x,t}$.

En esta secci\'on probaremos la existencia local de soluci\'on cl\'asica para la ecuaci\'on \eqref{Burg} con $\gamma\geq0$ y $\alpha\in(0,2]$ si el dato inicial est\'a en un cierto espacio de funciones. Para ello utilizaremos el m\'etodo de la energ\'ia de Leray (ver \cite{L}). El argumento de Leray ya es un argumento cl\'asico. Nosotros trataremos de aplicarlo con todo detalle en el ejemplo sencillo de la ecuaci\'on \eqref{Burg}.

La idea del m\'etodo es conseguir una sucesi\'on de problemas \emph{regularizados} para la ecuaci\'on \eqref{Burg}. Para todos los problemas regularizados se demuestra la existencia utilizando el Teorema de Picard en un espacio de Banach adecuado. Se concluye el argumento observando que la familia de soluciones regularizadas forma una sucesi\'on de Cauchy, y por lo tanto convergente en alg\'un espacio. El espacio de Banach que vamos a usar es $H^s(\RR)$, $s\geq3$ porque utilizando la inmersi\'on de Sobolev tenemos que entonces $u\in C^2(\RR)$. 

\subsection{Estimaciones \emph{a priori}}
Para utilizar el m\'etodo de Leray hemos de conseguir unas cotas \emph{a priori} para ciertas cantidades. Es decir, suponiendo que hay soluci\'on. 

Si pensamos en $u$ como en la velocidad de un fluido incompresible nos interesa conocer qu\'e ocurre con la \emph{'energ\'ia cin\'etica'} del fluido, es decir, la norma $L^2$ de $u$.

\begin{prop}[Principio del m\'aximo para $||\cdot||_{L^2}$]
\label{L2}
Sea $u$ una soluci\'on cl\'asica de \eqref{Burg}. Entonces se tiene que:
\begin{itemize}
\item[a)] Si $\gamma=0$,
$$|| u ||_{L^2} (t) = || f ||_{L^2},$$
\item[b)] Si $\gamma>0$, 
$$|| u ||_{L^2} (t) \leq ||f ||_{L^2}.$$
\end{itemize} 
\end{prop}
\begin{proof}
Observamos que 
\begin{multline*}
\frac{1}{2}\frac{d}{dt}||u||^2_{L^2}=\int_{\RR}u\pat udx=-\int_{\RR} u^2\pax udx-\gamma\int_{\RR}u\Lambda^{\alpha}udx\\
=-\frac{1}{3}\int_{\RR}\pax(u^3)dx-\gamma\int_{\RR}\Lambda^{\alpha}uudx=-\gamma\int_{\RR}\Lambda^{\alpha}uudx. 
\end{multline*} 
De aqu\'i se concluye la parte a) del resultado. Para obtener la parte b) hemos de utilizar el Teorema de Plancherel:
$$
\int_{\RR}\Lambda^{\alpha}uudx=\int_\RR|\xi|^{\alpha}\hat{u}\bar{\hat{u}}d\xi=\int_{\RR}\left(\Lambda^{\alpha/2}u\right)^2dx.
$$
De esta \'ultima igualdad (que no es m\'as que la prueba de que el operador $\Lambda^{\alpha}$ es sim\'etrico) se concluye el resultado.
\end{proof}

\begin{figure}[h]
	\centering
		\includegraphics[scale=0.30]{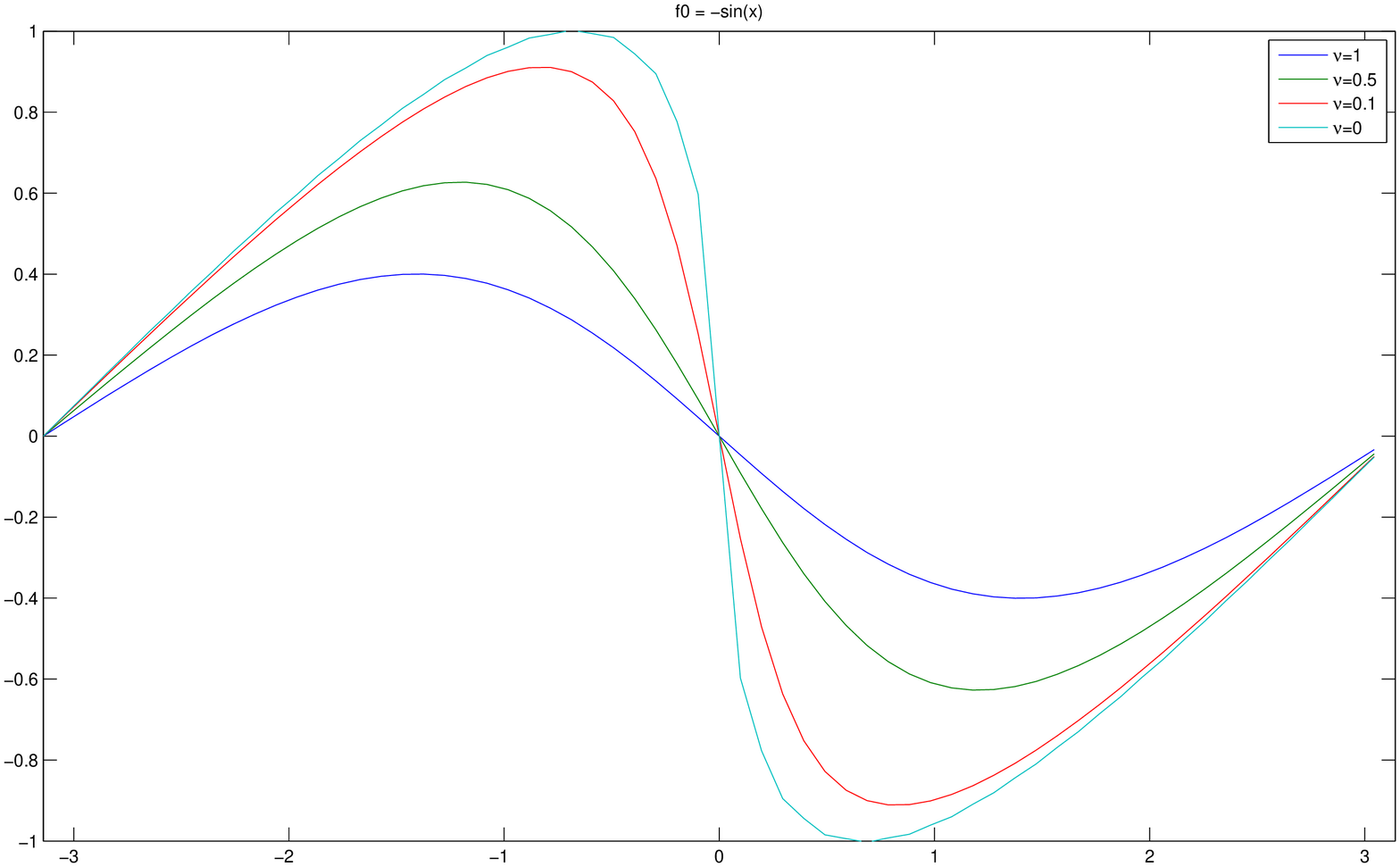} 
		\caption{Soluciones para distinas $\gamma$.}
\label{figura4}
\end{figure}

En el caso de que la soluci\'on $u$ no sea lo bastante regular como para ser soluci\'on cl\'asica entonces la norma $L^2$ no se conserva. Esto en las ecuaciones de Euler es una serie de celebrados art\'iculos (ver \cite{dLS}). Las soluciones d\'ebiles de \eqref{Euler} en tres dimensiones o de \eqref{Burg} con $\gamma=0$ no siempre tienen sentido f\'isico. Lo que implica f\'isicamente es que un fluido perfecto que estuviese inicialmente en reposo puede comenzar a agitarse a lo loco sin haber mediado fuerza externa alguna, siempre y cuyo la velocidad $u$ no tenga la suficiente regularidad como para ser soluci\'on cl\'asica. Matem\'aticamente esta \emph{'paradoja'} f\'isica se traduce en que no hay unicidad de soluciones d\'ebiles para el sistema de ecuaciones \eqref{Euler} en tres dimensiones. En el caso de \eqref{Burg} sin viscosidad ($\gamma=0$) lo que ocurre es que se da un choque de curvas caracter\'isticas y despu\'es no hay una manera \'unica de continuar la soluci\'on. Enlaza esto con las condiciones de \emph{Rankine-Hugoniot} y las \emph{soluciones de entrop\'ia} (ver \cite{Ev-08}). Adem\'as, si consideramos el problema viscoso 
$$
\pat u^\gamma+ u^\gamma\pax u^\gamma=-\gamma\Lambda^\alpha u^\gamma,
$$
pero tomamos $\gamma\rightarrow0$ recuperamos la soluci\'on de entrop\'ia para la ecuaci\'on de Burgers no viscosa como l\'imite, es decir, $u^0=\lim_{\gamma\rightarrow0} u^\gamma$ (ver Figura \ref{figura4}).

Ya tenemos una cota \emph{a priori} para la norma $L^2$. Como mencionamos antes queremos que nuestra soluci\'on est\'e en el espacio de Sobolev $H^3$, por lo que s\'olo falta estimar la norma $L^2$ de $\pax^3 u$.

\begin{prop} 
Sea $u$ una soluci\'on cl\'asica del problema \eqref{Burg} con un dato inicial $f\in H^3(\RR)$. Entonces se tiene la siguiente cota
\begin{equation}
||u||_{H^3}(t)\leq \frac{||f||_{H^3}}{1-c||f||_{H^3}t}.
\label{h3} 
\end{equation}
\end{prop}
\begin{proof}
Se tiene que 
$$
\frac{1}{2}\frac{d}{dt}||\pax^3 u||^2_{L^2}=\int_{\RR} \pax^3 u\pax ^3\pat udx,
$$ 
por lo que hemos de derivar tres veces la ecuaci\'on \eqref{Burg}, obteniendo
$$
 \pax^3\pat u=-\gamma\Lambda^{\alpha}\pax^3 u-3(\pax^2 u)^2+4\pax u\pax^3u+u\pax^4u.
$$
Introducimos esta expresi\'on en nuestra expresi\'on y logramos
\begin{eqnarray*}
\frac{1}{2}\frac{d}{dt}||\pax^3 u||^2_{L^2}&=&-\gamma\int_{\RR}\pax^3 u\Lambda^{\alpha}\pax^3u dx-3\int_{\RR}(\pax^2 u)^2\pax^3udx+\int_{\RR}4\pax u(\pax^3u)^2dx\\
&&+\int_\RR u\pax^4u\pax^3udx=I_1+I_2+I_3+I_4.\\
\end{eqnarray*}
Usando el Teorema de Plancherel se obtiene que 
$$
I_1=-\gamma\int_\RR|\xi|^{\alpha}\hat{\pax^3u}\bar{\hat{\pax^3u}}d\xi=-\gamma\int_{\RR}\left(\Lambda^{\alpha/2}\pax^3u\right)^2dx=-\gamma||\Lambda^{\alpha/2}\pax^3 u||_{L^2}^2\leq 0.
$$
La segunda integral se anula por las condiciones de borde impuestas
$$
I_2=-\int_{\RR}\pax\left[(\pax^2 u)^3\right]=0.
$$
Las integrales que faltan son las m\'as singulares por tener el mayor n\'umero de derivadas, sin embargo s\'olo nos hemos de preocupar de $I_3$, porque $I_4$ es igual. En efecto
$$
I_4=\frac{1}{2}\int_\RR u\pax\left[(\pax^3 u)^2\right]=-\frac{1}{2}\int_\RR \pax u(\pax^3 u)^2.
$$
Por lo tanto s\'olo hemos de acotar la integral $I_3$. Se tiene que 
$$
I_3\leq c||\pax u||_{L^\infty}||\pax^3 u||_{L^2}^2\leq ||u||_{H^3}^3,
$$
y conclu\'imos que 
$$
\frac{d}{dt}||\pax^3 u||^2_{L^2}\leq c||u||_{H^3}^3.
$$
Usando la Proposici\'on \ref{L2} obtenemos 
$$
\frac{d}{dt}||u||_{H^3}\leq c||u||_{H^3}^2.
$$
Sin m\'as que integrar la EDO se concluye la primera parte del resultado. 
\end{proof}
\noindent\begin{coment}
Es m\'as, si $f\in H^k(\RR)$, $k>3$ se tiene la siguiente cota
\begin{equation}
||u||_{H^k}(t)\leq \frac{||f||_{H^k}}{1-c||f||_{H^k}t}.
\label{hk} 
\end{equation}
\label{estim}
\end{coment}

\subsection{Regularizaci\'on del problema}
Consideremos una funci\'on $\rho\in C^\infty_c$, $\rho(x)=\rho(|x|)$, $\rho\geq0$ y $\int_{\RR}\rho=1$. Un ejemplo de tal $\rho$ puede ser la funci\'on 
$$
\rho(x)=\frac{e^{-\frac{1}{1-x^2}}}{\int_{-1}^1e^{-\frac{1}{1-y^2}}dy}.
$$
Ahora, para cualquier $\epsilon$ positivo, consideramos las funciones 
\begin{equation}
\rho_\epsilon=\frac{1}{\epsilon}\rho\left(\frac{x}{\epsilon}\right).
\label{epsi} 
\end{equation}
Observamos que $\reps$ siguen teniendo las mismas propiedades que $\rho$. Dada una funci\'on $g\in L^p$, si ahora hacemos la convoluci\'on $g_\epsilon=\rho_\epsilon *g$, obtenemos que $g_\epsilon\in C^\infty$ y, gracias a la desigualdad de Young:
\begin{equation}
||\reps* g||_{L^r}\leq ||\rho_\epsilon||_{L^q}||g||_{L^p},\quad \frac{1}{r}+1=\frac{1}{p}+\frac{1}{q},
\label{Young} 
\end{equation}
tenemos que 
$$
||\reps * g||_{L^p}\leq||g||_{L^p}.
$$
Adem\'as $\lim_{\epsilon\rightarrow0} g_\epsilon=g$. Pueden verse m\'as propiedades de estos n\'ucleos $\rho_\epsilon$ en \cite{bertozzi-Majda}.

El primer paso del m\'etodo es regularizar el problema de manera que podamos demostar existencia de los problemas regularizados y tambi\'en podamos usar las estimaciones \emph{a priori}. Para ello utilizarmos la convoluci\'on con $\rho_\epsilon$. As\'i consideramos la familia de problemas
\begin{equation}
\left\{\begin{array}{l}
\pat u_\epsilon=- \reps *\left(\reps *u_\epsilon\right)\pax \left(\reps *u_\epsilon\right)-\reps *\left(\gamma\Lambda^\alpha \reps * u_\epsilon\right),\;\; (t,x)\in[0,T]\times\RR, \epsilon>0\\
u_\epsilon(0,x)=f(x).\\
\end{array}\right.
\label{Burgeps}
\end{equation}
 
Esta regularizaci\'on particular se ha elegido pensando en obtener las mismas estimaciones de la Proposici\'on \ref{estim}. Por ser $\reps$ una funci\'on radial tenemos que 
$$
\int f\reps *g=\int\int f(x)\reps(x-y)g(y)dydx=\int\int f(x)\reps(y-x)g(y)dydx=\int g\reps * f,
$$
e integrando por partes tenemos unas estimaciones para $u_\epsilon$ como las de la Proposici\'on \ref{estim}.

\begin{prop} 
Sea $u_\epsilon$ soluci\'on de \eqref{Burgeps} con dato inicial $f\in H^k(\RR)$, $k\geq3$. Entonces se tiene la siguiente cota 
\begin{equation}
||u_\epsilon||_{H^k}(t)\leq \frac{||f||_{H^k}}{1-c||f||_{H^k}t}.
\label{hkeps} 
\end{equation}
\label{estimeps}
\end{prop}

\subsection{Existencia para los problemas regularizados}
Para demostrar la existencia de los problemas regularizados \eqref{Burgeps} vamos a utilizar el teorema de Picard en el espacio $H^k(\RR)$. Hemos de ver entonces que para todo $\epsilon$ se tiene que 
$$
F_\epsilon=- \reps *\left(\reps *u_\epsilon\right)\pax \left(\reps *u_\epsilon\right)-\reps *\left(\gamma\Lambda^\alpha \reps * u_\epsilon\right)
$$
es localmente (es decir, si $||u||_{H^k}<\lambda$ para cierto $\lambda$) Lipschitz con respecto a la norma $H^k$. Para empezar hemos de asegurarnos que no perdemos derivadas al aplicar $F_\epsilon$. Esto se consigue porque en lugar de derivar $u_\epsilon$ derivamos el n\'ucleo $\reps$, de manera que $u_\epsilon$ no pierde derivadas. Veamos que $F_\epsilon$ es Lipschitz:
\begin{lem} Sea $\lambda>0$, entonces para todo $\epsilon>0$ $F_\epsilon$ es una funci\'on Lipschitz en $\{g: g\in H^k(\RR), ||g||_{H^k}\leq \lambda\}$ con $k\geq 3$. 
\end{lem}
\begin{proof}Se tiene que 
\begin{multline*}
||F_\epsilon (u)-F_\epsilon(v)||_{H^k}= ||- \reps *\left(\reps *u\right)\pax \left(\reps *u\right)-\reps *\left(\gamma\Lambda^\alpha \reps * u\right)\\
+ \reps *\left(\reps *v\right)\pax \left(\reps *v\right)+\reps *\left(\gamma\Lambda^\alpha \reps * v\right)||_{H^k}.
\end{multline*}
Vamos a agrupar las partes difusivas de $u$ con la de $v$:
$$
||\reps *\left(-\gamma\Lambda^\alpha \reps * u+\gamma\Lambda^\alpha \reps * v\right)||_{H^k}\leq\gamma||\reps *(v-u)||_{H^{k+\alpha}}\leq L_1(\gamma,\alpha,\epsilon)||u-v||_{H^k}.
$$ 
Queda probar la estimaci\'on para la parte convectiva:
\begin{multline*}
|| \reps *(-\reps *u\pax \left(\reps *u\right) + \reps *v\pax (\reps *v))||_{H^k}\\
\leq|| -\reps *u\pax \left(\reps *u\right)\pm \reps *u\pax \left(\reps *v\right) + \reps *v\pax (\reps *v)||_{H^k}\\
\leq|| \reps *u\pax \left(\reps * (v-u)\right)+ (\reps *(v-u))\pax (\reps *v)||_{H^k}\\
\leq|| \reps *u\pax \left(\reps * (v-u)\right)||_{H^3}+ ||(\reps *(v-u))\pax (\reps *v)||_{H^k}.
\end{multline*}
Utilizamos ahora que, si $s>1/2$, $H^s(\RR)$ es un \'algebra de Banach, es decir, que se cumple
$$
||fg||_{H^s}\leq||f||_{H^s}||g||_{H^s}.
$$
Gracias a esta propiedad de los espacios de Sobolev tenemos que 
\begin{multline*}
|| \reps *u\pax \left(\reps * (v-u)\right)||_{H^k}+ ||(\reps *(v-u))\pax (\reps *v)||_{H^k}\\
\leq ||\reps *u||_{H^k}||\pax \left(\reps * (v-u)\right)||_{H^k}+||\pax (\reps *v)||_{H^k}||(\reps *(v-u))||_{H^k}\\
\leq L_2(\epsilon,\lambda)||u-v||_{H^k}.
\end{multline*}
Para concluir hemos de elegir $L=\max\{L_1,L_2\}$.
\end{proof}

Por lo tanto, dado un dato inicial $f\in H^k$, si aplicamos el Teorema de Picard tenemos que existe una sucesi\'on $u_\epsilon\in C^1([0,T_\epsilon],H^k(\RR))$. Adem\'as, por las propiedades de los n\'ucleos $\rho_\epsilon$, se puede demostrar que $T_\epsilon=\infty.$

Tenemos as\'i el siguiente resultado:
\begin{lem}
Dado una dato inicial $f\in H^k$ con $k\geq 3$ el problema \eqref{Burgeps} tiene una \'unica soluci\'on $u_\epsilon\in C^1([0,\infty),H^k)$. 
\end{lem}

\subsection{Paso al l\'imite $\epsilon\rightarrow0$}
Veremos que $u_\epsilon\in C^1([0,T],H^k)$ forma una sucesi\'on de Cauchy en el espacio $C([0,T],H^s)$, $0\leq s<k$ donde $T<T^*=\frac{1}{c||f||_{H^k}}$ (ver Proposici\'on \ref{estimeps}). Para ello tenemos que ver que en $L^2$ es Cauchy y entonces utilizando las estimaciones de la Proposici\'on \ref{estimeps} concluiremos que es una sucesi\'on de Cauchy en $H^s$ con $s<k$. El caso extremo de $s=k$ lo trataremos al final.

\begin{lem}
$u_\epsilon$ es una sucesi\'on de Cauchy en $C([0,T],H^s)$, con $0\leq s<k$. 
\end{lem}
\begin{proof}
Veremos primero que $u_\epsilon$ es Cauchy en $C([0,T],L^2)$. Sean $\epsilon, \delta$ dos n\'umeros positivos. Entonces se tiene
\begin{eqnarray*}
\frac{1}{2}\frac{d}{dt}||u_\epsilon-u_\delta||_{L^2}^2&=&-\gamma\int_{\RR}(\rho_\epsilon *(\Lambda^\alpha \reps * u_\epsilon)-\rho_\delta *(\Lambda^\alpha \rho_\delta * u_\delta))(u_\epsilon-u_\delta)dx\\
&&-\int_\RR(\reps *(\reps *u_\epsilon(\pax(\reps *u_\epsilon)))-\rho_\delta *(\rho_\delta * u_\delta(\pax(\rho_\delta * u_\delta))))(u_\epsilon-u_\delta)dx\\
&=&I_1+I_2.
\end{eqnarray*}
Veamos el caso de los operadores difusivos:
\begin{multline*}
I_1=-\gamma\int_{\RR}(\rho_\epsilon *(\Lambda^\alpha \reps * u_\epsilon)-\rho_\delta *(\Lambda^\alpha \rho_\delta * u_\delta))(u_\epsilon-u_\delta)dx\\
=-\gamma\int_{\RR}(\rho_\epsilon *(\Lambda^\alpha \reps * u_\epsilon)\pm\rho_\delta *(\Lambda^\alpha \rho_\delta * u_\epsilon)-\rho_\delta *(\Lambda^\alpha \rho_\delta * u_\delta))(u_\epsilon-u_\delta)dx\\
= -\gamma\int_{\RR}((\rho_\epsilon *(\Lambda^\alpha \reps) -\rho_\delta *(\Lambda^\alpha \rho_\delta))* u_\epsilon+\rho_\delta *(\Lambda^\alpha \rho_\delta * (u_\epsilon-u_\delta))(u_\epsilon-u_\delta)dx.
\end{multline*}
Integrando por partes el segundo sumando obtenemos que
\begin{multline*}
-\gamma\int_{\RR}((\rho_\epsilon *(\Lambda^\alpha \reps) -\rho_\delta *(\Lambda^\alpha \rho_\delta))* u_\epsilon+\rho_\delta *(\Lambda^\alpha \rho_\delta * (u_\epsilon-u_\delta))(u_\epsilon-u_\delta)dx\\
\leq-\gamma\int_{\RR}(\rho_\epsilon *(\Lambda^\alpha \reps) -\rho_\delta *(\Lambda^\alpha \rho_\delta))* u_\epsilon(u_\epsilon-u_\delta)\\
\leq \gamma ||u_\epsilon-u_\delta||_{L^2}||(\rho_\epsilon *(\Lambda^\alpha \reps) -\rho_\delta *(\Lambda^\alpha \rho_\delta))* u_\epsilon||_{L^2},
\end{multline*}
y por lo tanto, usando las propiedades de las suavizaciones (ver \cite{bertozzi-Majda})
\begin{eqnarray*}
||\Lambda^\alpha((\rho_\epsilon *( \reps) -\rho_\delta *( \rho_\delta))* u_\epsilon)||_{L^2}&\leq& ||(\rho_\epsilon *( \reps) -\rho_\delta *( \rho_\delta))* u_\epsilon||_{H^\alpha}\\
&=&||(\rho_\epsilon *( \reps) -\rho_\delta *( \rho_\delta))* u_\epsilon\pm \reps * u_\epsilon\pm \rho_\delta * u_\epsilon||_{H^\alpha}\\
&\leq& ||\rho_\epsilon *( \reps *u_\epsilon)-\reps * u_\epsilon||_{H^\alpha}\\
&&+||\rho_\delta *( \rho_\delta *u_\epsilon)-\rho_\delta * u_\epsilon||_{H^\alpha}\\
&&+||\reps * u_\epsilon-\rho_\delta * u_\epsilon||_{H^\alpha}\\
&\leq& c\epsilon||\rho_\epsilon * u_\epsilon||_{H^{1+\alpha}}+c\delta||\rho_\delta * u_\epsilon||_{H^{1+\alpha}}\\
&&+||\reps * u_\epsilon\pm u_\epsilon-\rho_\delta * u_\epsilon||_{H^\alpha}\\
&\leq& c\epsilon||\rho_\epsilon * u_\epsilon||_{H^{1+\alpha}}+c\delta||\rho_\delta * u_\epsilon||_{H^{1+\alpha}}\\
&&+(c\epsilon+c\delta)||u_\epsilon||_{H^{1+\alpha}}\\
&\leq& c\max\{\epsilon,\delta\}||u_\epsilon||_{H^3}
\end{eqnarray*}
Recordemos que, fijo $T<T^*$, tenemos la cota uniforme $||u_\epsilon||_{H^3}\leq C(T,||f||_{H^3})$ y por lo tanto tenemos que 
$$
I_1\leq C(T,||f||_{H^3},\gamma) ||u_\epsilon-u_\delta||_{L^2}\max\{\epsilon,\delta\}.
$$
Para la integral $I_2$ se hace igual. Hay que sumar y restar los t\'erminos 
$$
\rho_\delta *(\reps * u_\epsilon\pax (\reps * u_\epsilon)),\;\;\rho_\delta *(\rho_\delta * u_\epsilon\pax (\reps * u_\epsilon)),\;\;\rho_\delta *(\rho_\delta * u_\delta\pax (\reps * u_\epsilon))...
$$
y aplicar las mismas desigualdades (Sobolev, Cauchy-Schwartz...).

Una vez que hemos probado que $u_\epsilon$ es una sucesi\'on de Cauchy en $C([0,T],L^2)$, utilizando la cota uniforme en $\epsilon$ (Proposici\'on \ref{estimeps}) 
$$
||u_\epsilon||_{H^k}\leq C(T,||f||_{H^k}),
$$
podemos interpolar (ver \cite{bertozzi-Majda}) entre el espacio $H^0=L^2$ y $H^k$, de manera que $u_\epsilon$ es de Cauchy para todo $H^s$ con $s<k$. Podemos as\'i asegurar la existencia de $u\in C([0,T],H^s)$ como el l\'imite de $u_\epsilon$ en la topolog\'ia usual de ese espacio.
\end{proof}

Adem\'as podemos mejorar el resultado para cubrir el caso extremo $s=k$. As\'i hemos de ver que $u=\lim_{\epsilon\rightarrow 0}u_\epsilon$ est\'a en el espacio $H^k$.\footnote{Podemos elegir la noci\'on de l\'imite de cualquier espacio $C([0,T],H^s)$ con $s<k$.} 

\begin{lem}
La funci\'on $u$ definida como el l\'imite anterior est\'a en $C([0,T],H^k)$. 
\end{lem}
\begin{proof}
Para probarlo hemos de razonar con an\'alisis funcional. Como $\{u_\epsilon\}$ est\'a acotada en $H^k$ (Proposici\'on \ref{estimeps}) tiene una subsucesi\'on que converge d\'ebilmente en este espacio. Es decir
$$
<u_\epsilon,z>\rightarrow <v,z>,\;\;\forall z\in H^{-k}.
$$
En principio no podemos asegurar que $u=v$. Para concluir la igualdad utilizamos las convergencias d\'ebiles, pues la convergencia fuerte implica la convergencia d\'ebil, en los espacios $H^s$ con $s<k$ adem\'as del hecho de que $H^k\subset H^s\Rightarrow H^{-k}\supset H^{-s}$ y obtenemos que si $u\neq v$ tenemos una contradicci\'on.
\end{proof}

\subsection{Regularidad de $u$}
Gracias a los lemas anteriores podemos demostrar el siguiente resultado:
\begin{teo}[Existencia local de soluci\'on cl\'asica]
Sea $f\in H^k$, $k\geq3$. Entonces existe una \'unica $u\in C([0,T],H^k)\cap C^1([0,T],H^{k-2})$, definida como el l\'imite uniforme en compactos de $\RR$ de $u_\epsilon$, soluci\'on cl\'asica del problema \eqref{Burg} con $\gamma\geq0$ y $0<\alpha\leq 2$.
\end{teo}
\begin{proof}

\emph{Existencia:}
De los lemas anteriores se obtiene la existencia de $u\in C([0,T],H^k)$ como l\'imite uniforme en compactos de $u_\epsilon$. Para obtener la regularidad temporal necesaria para $u$ observamos que $\pat u_\epsilon$ tiende como distribuci\'on a la derivada d\'ebil de $u$ (que denotamos por $\pat u$). En efecto, para toda funci\'on test $\phi$, se tiene que 
$$
\int_\RR\int_0^T\phi\pat u_\epsilon dtdx=-\int_\RR\int_0^T\pat \phi u_\epsilon dt dx\rightarrow -\int_\RR\int_0^T\pat \phi udtdx= \int_\RR\int_0^T\phi \pat u dt dx.
$$
Adem\'as $\pat u_\epsilon\rightarrow -\gamma\Lambda^{\alpha} u- u \pax u$ en $C([0,T], C(\RR))$ (y por lo tanto tambi\'en en sentido distribucional). Esto es una consecuencia de las propiedades de los n\'ucleos suavizantes.

Ahora estamos en condiciones de conseguir una mejor estimaci\'on para $\pat u$:
\begin{multline*}
||\pat u||_{C([0,T], C(\RR))}=\max_{0\leq t\leq T} \max _{x\in U;\newline U\subset\subset \RR}|-\gamma\Lambda^\alpha u- u\pax u|\leq \max_{0\leq t\leq T}\left(c||u||_{C^2}+ ||u||_{C^2}^2\right)\\
\leq c||u||_{C([0,T], H^3(\RR))}+ ||u||_{C([0,T], H^3(\RR))}^2\leq C(||f||_{H^k},T,\gamma).
\end{multline*}
En particular si $u\in C([0,T],H^k)$ con $k\geq 3$ se concluye que $u$ es soluci\'on cl\'asica de la ecuaci\'on de Burgers.

Para concluir que $u\in C^1([0,T],H^{k-2})$ tenemos que obtener una cota para $||\pat u||_{H^{k-2}}$: 
$$
||\pat u||_{H^{k-2}}\leq \gamma||\Lambda^\alpha u||_{H^{k-2}}+||u||_{H^{k-2}}||\pax u||_{H^{k-2}}\leq C(||f||_{H^k},T,\gamma),
$$
y entonces
\begin{multline*}
||u||_{C^1([0,T],H^{k-2})}=\max_{0\leq t\leq T}||u||_{H^{k-2}}+||\pat u||_{H^{k-2}}\\
\leq \max_{0\leq t\leq T}||u||_{H^{k-2}}+ \max_{0\leq t\leq T}\gamma||\Lambda^\alpha u||_{H^{k-2}}+\max_{0\leq t\leq T}||u||_{H^{k-2}}||\pax u||_{H^{k-2}}\leq C(||f||_{H^k},T,\gamma),
\end{multline*}
donde en la \'ultima desigualdad hemos usado \eqref{hk}.

\emph{Unicidad:}
Supongamos que hubiese dos soluciones, $u$ y $v$ en el espacio $H^k$ con $k\geq3$ del problema \eqref{Burg} con el mismo dato inicial. Entonces 
\begin{eqnarray*}
\frac{1}{2}\frac{d}{dt}||u-v||_{L^2}^2&=&-\gamma||\Lambda^{\alpha/ 2}(u-v)||_{L^2}^2-\int_\RR (u-v)(u\pax u-v\pax v)dx\\
&\leq&-\frac{1}{2}\int_\RR (u-v)\pax(u^2-v^2)dx\\
&\leq&\frac{-1}{2}\int_\RR(u-v)\pax((u+v)(u-v))dx\\
&\leq&\frac{-1}{2}\int_\RR(u-v)^2\pax(u+v)dx-\frac{1}{2}\int_\RR(u-v)\pax(u-v)(u+v)dx\\
&\leq& c\int_\RR (u-v)^2\pax (u+v)dx\\
&\leq& c(||u||_{H^k},||v||_{H^k})||u-v||_{L^2}^2\\
&\leq& c(T,||f||_{H^k})||u-v||_{L^2}^2,
\end{eqnarray*}
y usando Gronwall se obtiene la unicidad.
\end{proof}

\section{Propiedades de las soluciones}
En esta secci\'on encontraremos algunas propiedades b\'asicas que deber\'a tener una soluci\'on de la ecuaci\'on de Burgers.   

Como la ecuaci\'on \eqref{Burg} modeliza un escalar que es transportado por un fluido donde la velocidad del fluido la hemos hecho proporcional al escalar \eqref{Burgers} es de esperar que la \emph{masa total}, \emph{i.e.} se conserva. En efecto,
\begin{prop}[Conservaci\'on de la masa]
Sea $u$ una soluci\'on cl\'asica de la ecuaci\'on \eqref{Burg} entonces se tiene que 
$$
\int_\RR u(t,x)dx=\int_\RR f(x)dx.
$$ 
\end{prop}
\begin{proof}
Integrando la ecuaci\'on \eqref{Burg} en espacio 
$$
\frac{d}{dt}\int_\RR u(t,x)dx=\int \pat u(t,x)dx=\int_\RR-\frac{1}{2}\pax(u(t,x)^2)-\gamma\Lambda^\alpha u(t,x)dx=I_1+I_2.
$$ 
Por las condiciones de borde para $u$ la integral $I_1$ se anula.

Para ver que la integral $I_2$ se anula razonamos utilizando la transformada de Fourier. Dada un funci\'on $g(x)$ se cumple que 
$$
\int_\RR g(x)dx=\hat{g}(0).
$$
Como $\Lambda^\alpha$ en el espacio de Fourier es el multiplicador de la definici\'on \eqref{Lambdaalpha} obtenemos que 
$$
I_2=-\gamma\int_\RR \Lambda^\alpha u(t,x)dx=-\gamma\widehat{\Lambda^{\alpha}u}(0)=0.
$$ 
\end{proof}

Otra cantidad de la que es interesante ver la evoluci\'on es la norma $L^\infty$ de la soluci\'on:

\begin{prop}[Principio del m\'aximo para $||\cdot||_{L^\infty}$]\label{prop2}
Sea $u$ una solucion cl\'asica de \eqref{Burg} tal que $\max_{x\in\RR}u(t,x)\geq0$, $\min_{x\in\RR}u(t,x)\leq0$ y $\alpha\in (0,2]$, entonces se tienen los siguientes resultados
\begin{itemize}
\item[a)] Si $\gamma=0$,
$$|| u ||_{L^\infty} (t) = || f ||_{L^\infty},$$
\item[b)] Si $\gamma>0$, 
$$|| u ||_{L^\infty} (t) \leq ||f ||_{L^\infty}.$$
\end{itemize}
\end{prop}
\begin{proof}
Observamos que la funci\'on $M(t)=\max_{x\in\RR}u(t,x)=u(t,x_t)$ es Lipschitz. En efecto:
$$
\max_x |u(t_1,x)|=\max_x (|u(t_1,x)-u(t_2,x)+u(t_2,x)|)\leq \max_x (|u(t_1,x)-u(t_2,x)|)+\max_x |u(t_2,x)|,
$$
de donde
\begin{multline}
|\max_x u(t_1,x)-\max_x u(t_2,x)|\leq \max_x (|u(t_1,x)-u(t_2,x)|)=\max_x(|\pat u (s,x)||t_1-t_2|)\\ \leq\max_{s\in (t_2,t_1)}\max_x(|\pat u (s,x)|)|t_1-t_2|.
\label{eq6}
\end{multline}
Para concluir que la funci\'on es Lipschitz observamos que debemos restringirnos a un intervalo temporal $[0,T]$ con $T>0$ fijo. Ahora nos aseguramos de que $(t_1,t_2)\subset[0,T]$ y obtenemos
$$
|\max_x u(t_1,x)-\max_x u(t_2,x)|\leq\max_{s\in [0,T]}\max_x(\pat u (s,x))(t_1-t_2)=L(t_1-t_2).
$$
\noindent Usando el Teorema de Radamacher tenemos que $M(t)$ es derivable en casi todo punto,
\begin{eqnarray*}M'(t) &=& \lim_{h_j\rightarrow 0} \frac{M(t+h_j) - M(t)}{h_j} = \lim_{h_j\rightarrow 0}\frac{u(x_{t+h_j},t+h_j) - u(x_t,t)}{h_j}\\
&=&\lim_{h_j\rightarrow 0}\frac{u(x_{t+h_j},t+h_j)\pm u(x_{t+h_j},t) - u(x_t,t)}{h_j}= \pax u(t,x_t)+\pat u(x_t,t)\\
&=&\pat u(x_t,t).
\end{eqnarray*}
De manera que, si $\gamma=0$,
$$M'(t) = \pat u(x_t, t) = -u(x_t,t) \pax u(x_t,t) = -\frac{1}{2} \pax(u^2(x_t,t)) = 0$$
donde en la \'ultima igualdad hemos usado que $u(x_t,t)$ es m\'aximo. 
En el caso en el que $\gamma\geq0$ tenemos que dar signo al t\'ermino $-\gamma\Lambda^\alpha u(t,x_t)$. Usaremos la expresi\'on como una convoluci\'on \eqref{conv}. Se tiene que 
$$
-\gamma \Lambda^{\alpha}u(x_t,t) = c\text{P.V.}\int_{\RR}\frac{u(y,t) - u(x_t,t)}{|x-y|^{1+\alpha}}dy \leq 0.
$$
En el caso en el que $\alpha=2$ el resultado se reduce a conocer qu\'e signo tiene la segunda derivada en un punto de m\'aximo. 

Para el m\'inimo se hace de manera an\'aloga.
\end{proof}

\section{Un resultado an\'alogo al de Beale-Kato-Majda}
El lector atento se habr\'a percatado de que la ecuaci\'on de Burgers \eqref{Burg} tiene una especie de \emph{Teorema de Beale-Kato-Majda}. Recordemos que en una dimensi\'on espacial lo m\'as parecido a $\omega$ es la primera derivada espacial $\pax u$. Entonces tenemos el siguiente resultado:
\begin{prop}
Sea $T>0$. Entonces, si
$$
\int_0^T ||\pax u||_{L^\infty}(s)ds<\infty,
$$ 
existe $u$ soluci\'on cl\'asica de \eqref{Burg} al menos hasta tiempo $T$.
\end{prop}
\begin{proof}
En la prueba de la proposici\'on de las estimaciones \emph{a priori} (Proposici\'on \ref{estim}) logramos la cota 
$$
\frac{d}{dt}||\pax^3 u||_{L^2}^2\leq c||\pax u||_{L^\infty}||\pax^3 u||_{L^2}^2.
$$
De esta expresi\'on obtenemos 
$$
\frac{d}{dt}||u||_{H^3}\leq c||\pax u||_{L^\infty}|| u||_{H^3},
$$
y aplicando la desigualdad de Gronwall (es decir, integramos la EDO) conseguimos el resultado.
\end{proof}

\section{Buscando singularidades}
Consideramos la ecuaci\'on \eqref{Burg} con $\gamma=0$. Esta ecuaci\'on es una ecuaci\'on de transporte unidimensional no-lineal, por lo tanto es susceptible de aplicarle el m\'etodo de las caracter\'isticas. Esa t\'ecnica aplicada a esta ecuaci\'on se puede encontrar en casi todos los manuales de ecuaciones diferenciales. Nosotros utilizaremos una t\'ecnica similar a la que utilizamos para probar el Principio del M\'aximo en la secci\'on 2. Queremos encontrar una singularidad para $u$. Por la Proposici\'on \ref{prop2} tenemos que $u$ no \emph{explota}, por lo tanto si hay una singularidad debe estar en alguna de las derivadas espaciales de $u$. Como estamos en el caso no viscoso solamente tenemos una derivada espacial, así que $\pax u$ es nuestra candidata a cantidad que explota.

As\'i sea $u\in H^3(\RR)$ una soluci\'on cl\'asica de \eqref{Burg}. Por la inmersi\'on de Sobolev se tiene que $\pax u\in C^1(\RR)$. Ahora la idea es aplicar el Teorema de Rademacher a $\pax u$. De manera an\'aloga a la de la demostraci\'on de la Proposici\'on \ref{prop2} se demuestra que si $m_x(t)=\min_{x\in\RR}\pax u=\pax u (x_t,t)$ se tiene que 
$$
m'_x(t)=\pat\pax u(x_t,t)=-(\pax u(x_t,t))^2-u(x_t,t)\pax^2u(x_t,t).
$$
Como $x_t$ es el punto de m\'aximo de $\pax u$ tenemos que $\pax^2u=0$, por lo que la ecuaci\'on para la evoluci\'on del m\'aximo es
$$
m_x'=-m_x^2.
$$
Si ahora resolvemos la EDO obtenemos que 
\begin{equation}
m_x(t)=\frac{m_x(0)}{1+tm_x(0)}.
\label{sing} 
\end{equation}
Hemos obtenido el siguiente resultado 
\begin{prop}[Singularidad]
Sea $f\in H^3$ tal que $\min_{x\in \RR}\pax f(x)<0$. Entonces $u$ soluci\'on cl\'asica de \eqref{Burg} con $\gamma=0$ desarrolla una singularidad en tiempo finito.
\label{singula}
\end{prop}

\subsection{Interpretando la singularidad}
\begin{figure}[h]
	\centering
		\includegraphics[scale=0.4]{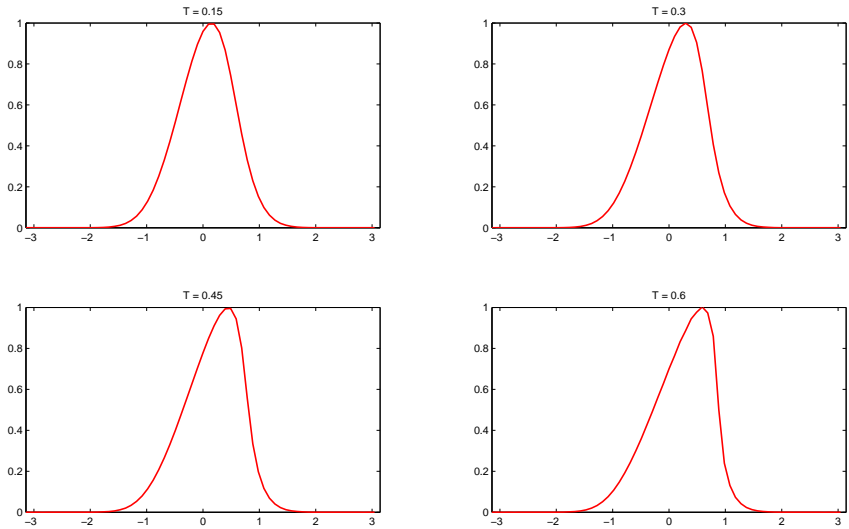}
		\caption{Simulaci\'on de una ola.}
\label{figura2}
\end{figure}

El tipo de singularidad que se produce se conoce como \emph{'choque'}\footnote{Del ingl\'es \emph{'shock'}.} y es una discontinuidad en $u$, o lo que es lo mismo $\pax u(x,T^*)=-\infty$. Si pensamos en $u$ como la altura de la superficie de una ola unidimensional lo que ocurre es que las part\'iculas del fluido que integran la ola se mueven m\'as r\'apido cuanto m\'as altas est\'en, por lo que las m\'as altas van a acercarse a las m\'as bajas. Entonces cuando las alcancen se produce este \emph{choque}. Esta discontinuidad previene que nuestra interfase $u$ deje de ser un grafo, porque si pensamos en la f\'isica de la ola tenemos que las part\'iculas m\'as r\'apidas sobrepasan a las m\'as lentas (ver Figura \ref{figura2}).

Hemos mencionado que la ecuaci\'on de Burgers pod\'ia ser entendida como una ola, lo que no es completamente descabellado. 

Si $y=u(x,t)$ es la superficie del agua (entonces la ola viene dada por la curva $(x,u(x,t))$), y $v=(v_1,v_2)$ es el campo de velocidades del agua, entonces la ecuaci\'on de la interfase entre el agua y el aire es
$$
\pat u=-\pax u v_1(x,u(x,t),t)+v_2(x,u(x,t),t)=(-\pax u,1)\cdot v.
$$
Ahora observamos que en las olas la parte m\'as alta se mueve a mayor velocidad que la parte m\'as baja. Por lo tanto parece natural (al menos como primera aproximaci\'on) hacer la hip\'otesis $v_1(x,u(x,t),t)=cu(x,t)$ (la constante $c$ podemos tomarla igual a 1 si hacemos un cambio de variables). As\'i si despreciamos la velocidad vertical, $v_2(x,u(x,t),t)\approx0$, tenemos la ecuaci\'on de Burgers no viscosa como modelo del perfil de una ola. Podemos razonar tambi\'en que la curva cambia seg\'un sea su curvatura, de manera que $v_2(x,u(x,t),t)\approx \kappa[u]$ con $\kappa[u]$ la curvatura de la curva $u$. Adem\'as si suponemos curvas suaves de amplitud peque\~ na podemos linealizar la curvatura obteniendo el operador $\gamma\pax^2 u(x,t)$. As\'i, con hip\'otesis razonables al menos en primera aproximaci\'on, obtenemos la ecuaci\'on de Burgers viscosa como modelo simplificado del perfil de una ola.

Una ola unidimensional con unas ciertas caracter\'isticas (longitud de onda, amplitud...) puede modelizarse en un mejor nivel de aproximaci\'on con la ecuaci\'on KdV siguiente:
\begin{equation}
\pat u+\pax uu=\pax^3u.
\label{KdV}
\end{equation}
Veamos de manera muy resumida c\'omo puede derivarse formalmente esta ecuaci\'on. Para describir una ola consideramos que el agua bajo la superfie tiene un flujo irrotacional, $i.e.$ $v=\nabla\phi$ para cierta funci\'on escalar $\phi$. Si suponemos v\'alidas las ecuaciones de Euler para el agua bajo la superficie tenemos que $\phi$ sigue la ley
$$
\pat \phi+\frac{1}{2}|\nabla\phi|^2+p+gy=0.
$$
Adem\'as, por la incompresibilidad se tiene $\Delta \phi=0$.


\begin{figure}
\center
\scalebox{1} 
{
\begin{pspicture}(0,-1.545)(6.24,1.565)
\psbezier[linewidth=0.04](0.0,0.735)(0.0,-0.065)(0.92163783,0.97220963)(1.92,0.915)(2.9183621,0.85779035)(3.1040077,0.5444405)(4.1,0.455)(5.095992,0.3655595)(6.176627,-0.48499432)(6.18,0.515)
\psline[linewidth=0.04cm](0.06,-1.485)(0.06,-1.505)
\psline[linewidth=0.04cm](0.08,-1.525)(6.22,-1.505)
\psline[linewidth=0.04cm](0.08,0.435)(6.1,0.435)
\psline[linewidth=0.04cm](1.24,-1.485)(1.24,0.355)
\usefont{T1}{ptm}{m}{n}
\rput(1.4478126,-0.9){h}
\psline[linewidth=0.04cm](1.86,0.535)(1.86,0.855)
\usefont{T1}{ptm}{m}{n}
\rput(2.013125,0.66){a}
\psline[linewidth=0.04cm](0.16,1.095)(5.8,1.115)
\usefont{T1}{ptm}{m}{n}
\rput(3.3092186,1.38){$\lambda$}
\end{pspicture}
}
\caption{Una ola y sus par\'ametros caracter\'isticos}
\end{figure}
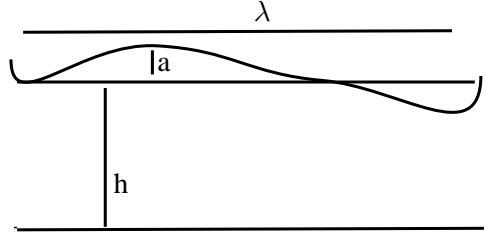

La coordenada $y$ se distingue de la $x$ en que act\'ua la gravedad, por lo tanto parece natural hacer un desarrollo de $\phi$ en potencias de $y$, $\phi=\sum_{n=0}^\infty y^n\phi_n(x,t)$. Si consideramos olas peque\~ nas en amplitud con respecto a la longitud de onda, entonces tenemos que despreciar los t\'erminos de orden grande en $y$ en nuestra expresi\'on para $\phi$.

Si adem\'as suponemos que $\pax\phi_0\approx u$ podemos concluir la ecuaci\'on \eqref{KdV}. Esta hip\'otesis se motiva por los desarrollos en serie de potencias anteriores.

Observamos que la hip\'otesis para obtener la ecuaci\'on de Korteveg-de Vries es menos restrictiva que para obtener la ecuaci\'on de Burgers, pues exclusivamente suponemos que $\pax\phi_0(x,t)=u(x,t)$ no que $\pax\phi(x,f(x,t),t)=v_1(x,f(x,t),t)=u(x,t)$. Consideramos por lo tanto discrepancias en los ordenes mayores.

\section{M\'etodo num\'erico}
En la siguiente secci\'on nos proponemos dar una aproximaci\'on num\'erica al problema \eqref{Burgers} propuesto en $\TT$.

\begin{figure}[h]
	\centering
		\includegraphics[scale=0.4]{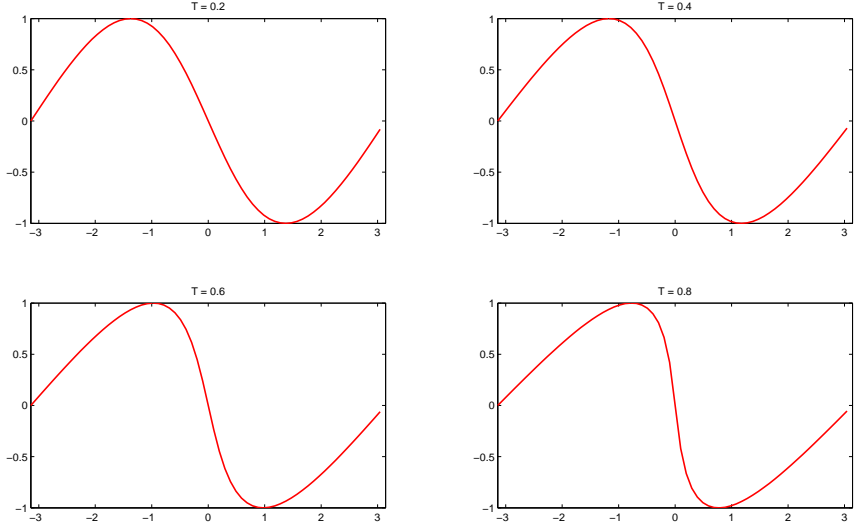} 
		\caption{Soluci\'on para $u(0,x)=-\sin(x)$ sin viscosidad.}
\label{figura3}
\end{figure}

\subsection{Discretizaci\'on espacial}
La forma natural de abordar el problema es usando la transformada de Fourier discreta (DFT) (ver \cite{can} y \cite{T}). Sea $\mathbb{T}^\star_N$ una discretizaci\'on de $\mathbb{T}$, con un n\'umero par $N$ de puntos, dada por,
$$x_j = \frac{\pi}{N}(2j-N),\quad j=0,\ldots,N-1$$
La idea es sustituir el valor de $u(\cdot,x) = u(x)$ en los nodos, $u(x_j)$, por su interpolante trigonom\'etrico,
\begin{eqnarray*}
I_N u(x) &=& \sum_{k={-N/2}}^{N/2-1} \tilde{u}_k e^{ikx}\\
I_N u(x_j) &=& u(x_j),\quad j=0,\ldots,N-1
\end{eqnarray*}
Es decir, dado un vector con los valores nodales de $u$, $U = [u_0\cdots u_{N-1}]$, la DFT nos proporciona $\tilde{U} = [\tilde{u}_{-N/2} \cdots \tilde{u}_{N/2-1}]$. Estos coeficientes vienen dados por,
$$\tilde{u}_k = \frac{1}{N} \sum_{j=0}^{N-1} u(x_j) e^{ikx_j},\quad j=0,\ldots,N-1$$
y puede demostrarse que est\'an relacionados con los coeficientes $\hat{u}_k$ de la serie de Fourier de $u(x)$ mediante,
$$\tilde{u}_k = \hat{u}_k + \sum_{m\neq0,m=-\infty}^{\infty} \hat{u}_{k+Nm},\quad k=-N/2,\ldots,N/2-1$$
Sabiendo como aproximar $u(x)$, es importante saber como aproximar $\pax u(x)$. Para ello derivamos el interpolante trigonom\'etrico anterior obteniendo,
$$\pax I_N u(x) = \sum_{k=-N/2}^{N/2-1} ik \tilde{u}_k e^{ikx}$$
Nuestra aproximaci\'on de $\pax u(x)$ en $x_l$ vendr\'a dada por lo que nosotros llamaremos,
$$(\mathcal{D}_N u)_l = \sum_{k=-N/2}^{N/2-1} ik\tilde{u}_k e^{-ikx_l},\quad l=0,\ldots,N-1$$
De la misma forma, parece apropiado aproximar el operador $\Lambda^\alpha = (-\pax^2)^{\alpha/2}$ mediante,
$$(\mathcal{E}_N^\alpha u)_l = \sum_{k=-N/2}^{N/2-1} |ik|^\alpha\tilde{u}_k e^{-ikx_l},\quad l=0,\ldots,N-1.$$

\begin{figure}[h]
	\centering
		\includegraphics[scale=0.4]{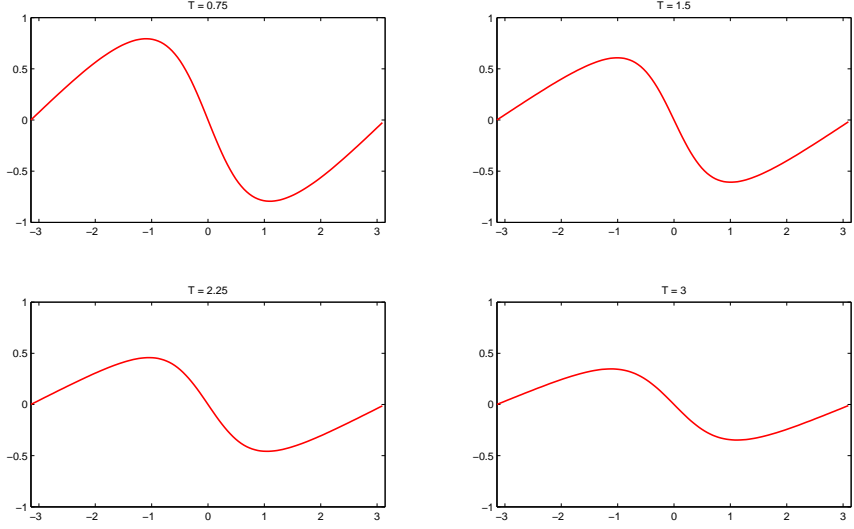}
		\caption{Soluci\'on para $u(0,x)=-\sin(x)$ con viscosidad.}
\label{figura5}
\end{figure}

\subsection{Discretizaci\'on temporal}
Para la aproximaci\'on de la derivada temporal en la ecuaci\'on utilizaremos el m\'etodo de Runge-Kutta expl\'icito de orden $4$, tambi\'en denominado RK4. 

Sea $T^\star = \{t_0,\ldots,t_M\}$ una discretizaci\'on del dominio temporal $[0,T]$. As\'i si $U_s= [u(x_1,t_s)\ldots u(x_N,t_s)]$ es la aproximaci\'on de la soluci\'on a nuestro problema en tiempo $t=t_s$ se tiene que $U_s$ verifica la ecuaci\'on
$$
\frac{d}{dt}U_s=F(U_s,t_s),
$$
con $F = - u\pax u - \gamma \Lambda^\alpha u$.

Entonces,
\begin{eqnarray*}
U_{s+1} = U_s + \frac16 \Delta t (K_1 + 2K_2 + 2K_3 + K_4),\\
\end{eqnarray*}
donde,
\begin{equation*}
\begin{tabular}{lll}
$K_1$& $=$ & $F(U_{s}, t_s)$,\\
$K_2$& $=$ & $F(U_{s} + \frac12 \Delta t K_1, t_s + \frac12 \Delta t)$,\\
$K_3$& $=$ & $F(U_{s} + \frac12 \Delta t K_2, t_s + \frac12 \Delta t)$,\\
$K_4$& $=$ & $F(U_{s} + \Delta t K_3, t_s + \Delta t)$.
\end{tabular}
\end{equation*}

\end{document}